\documentclass[11pt,leqno]{article}
\AtBeginDvi{}

\usepackage{amssymb}
\usepackage{amscd}
\usepackage{amsmath}
\usepackage{amsthm}
\usepackage{mathrsfs}
\usepackage{graphics}
\usepackage[dvipdfmx]{graphicx}
\usepackage{color}

\usepackage[all]{xy}
\usepackage{titlesec}
\titleformat{\section}[hang]
{\large\bf}
{\bf\thesection.}{.5em}{\bf }[ \vspace{+2pt}]
\titlespacing{\section}{0pt}{1.5ex plus .1ex minus .2ex}{0pt}
\titleformat{\subsection}[hang]
{\bf}
{\bf\thesubsection.}{.5em}{\bf }[ \vspace{+2pt}]
\titlespacing{\subsection}{0pt}{1ex plus .1ex minus .2ex}{0pt}
\def\upplim_#1{\underset{#1}{\overline\lim}\;}
\def\lowlim_#1{\underset{#1}{\underline\lim}\;}
\def\ang#1{{\langle}#1{\rangle}}
\makeatletter
\def\@makefnmark{\hbox{\@textsuperscript{\normalfont
\@thefnmark)}}}
\renewenvironment{enumerate}
  {\ifnum \@enumdepth >3\relax\@toodeep\else
   \advance\@enumdepth\@ne
   \edef\@enumctr{enum\romannumeral\the\@enumdepth}%
   \list{\csname label\@enumctr\endcsname}{%
         \ifnum \@listdepth=\@ne \topsep.1\normalbaselineskip
           \else\topsep\z@\fi
         \parskip\z@ \itemsep\z@ \parsep\z@
         \labelwidth1pc \labelsep0.5pc
         \ifnum \@enumdepth=\@ne \leftmargin1pc\relax
           \else\leftmargin\leftskip\fi
         \advance\leftmargin +1.5pc
         \usecounter{\@enumctr}%
         \def\makelabel##1{\hss\llap{##1}}}%
   \fi}{\endlist}
\makeatother
\newtheorem{cor}[equation]{Corollary}
\newtheorem{defn}[equation]{\indent{\it Definition}\rm }

\newtheorem{lem}[equation]{Lemma}
\newtheorem{llst}[equation]{\rm }
\newtheorem{prop}[equation]{Proposition}
\newtheorem{rmk}[equation]{\indent \rm {\it Remark}}
\newtheorem{thm}[equation]{Theorem}

\newcommand{\C}{{\mathbf{C}}}

\newcommand{\codim}{{\mathrm{codim}}}

\newcommand{\exh}{\hskip1pt{\widehat{\mathrm{ex}}}}

\newcommand{\del}{{\partial}}
\newcommand{\delbar}{\bar{\partial}}

\newcommand{\ep}{\varepsilon}
\newcommand{\fa}{^\forall}

\newcommand{\what}{\widehat}
\newcommand{\iso}{\cong}

\newcommand{\Jk}{{\what{J}_{k,A}}}
\newcommand{\Jkb}{\hskip3pt{\overline{\hskip-3pt\what{J}}_{k,A}}}
\newcommand{\Job}{\hskip3pt{\overline{\hskip-3pt\what{J}}_{0,A}}}
\newcommand{\sI}{\mathscr{I}}
\newcommand{\im}{{\Im\,}}
\newcommand{\Ker}{{\mathrm{Ker}\,}}

\newcommand{\Lie}{{\mathop{\mathrm{Lie}}}}
\newcommand{\lto}{\longrightarrow}

\newcommand{\N}{\mathbf{N}}

\renewcommand{\O}{{\mathcal{O}}}

\newcommand{\ord}{{\mathrm{ord}}}
\renewcommand{\P}{{\mathbf{P}}}

\newcommand{\vp}{\varphi}

\newcommand{\pnc}{{\mathbf{P}^n(\mathbf{C})}}

\newcommand{\ponec}{{\mathbf{P}^1(\mathbf{C})}}
\newcommand{\ptwoc}{{\mathbf{P}^2(\mathbf{C})}}

\newcommand{\Q}{{\mathbf{Q}}}

\newcommand{\re}{{\Re\,}}

\newcommand{\st}{\mathrm{St}}

\newcommand{\supp}{\mathrm{Supp}\,}
\newcommand{\td}{{\mathrm{tr.}\deg}\hskip0.3pt}
\newcommand{\tensor}{\otimes}

\newcommand{\zar}{\mathrm{Zar}}
\newenvironment{pf}{\par{\it Proof.\,}}{\qed\vskip+2pt}
\numberwithin{equation}{section}
\def\labelenumi{\rm(\roman{enumi})}

\title{Analytic Ax-Schanuel Theorem for semi-abelian varieties\\
and Nevanlinna theory}
\author{J. Noguchi\thanks{Research supported in part by Grant-in-Aid
 for Scientific Research (C) 19K03511, Japan Society for the Promotion of Science.}}
\date{\normalsize Graduate School of Mathematical Sciences\\
The University of Tokyo, Tokyo, Japan\\
e-mail: noguchi@g.ecc.u-tokyo.ac.jp}
\date{\it\today}
\begin{document}
\AtBeginDvi{}

\setlength{\baselineskip}{13pt}
\maketitle
\thispagestyle{empty}
\begin{abstract}
The purpose of this paper is to explore
Nevanlinna theory of the entire curve
 $\exh_A f:=(\exp_Af,f):\C \to A \times \Lie(A)$ associated with an
entire curve $f: \C \to \Lie(A)$,
 where $\exp_A:\Lie(A)\to A$ is an exponential map of
 a semi-abelian variety $A$.
Firstly  we  give a Nevanlinna theoretic proof
to the {\em analytic Ax-Schanuel Theorem} for semi-abelian varieties,
which was proved by J. Ax 1972 in the case of formal power series
 (Ax-Schanuel Theorem).
We assume some non-degeneracy condition for $f$ such that
the elements of the vector-valued function
 $f(z)-f(0) \in \Lie(A)\iso \C^n$
are $\Q$-linearly independent in the case of $A=(\C^*)^n$.
Then by making use of the Log Bloch-Ochiai Theorem and a key estimate
which we show, we  prove that $\td_\C\, \exh_A f \geq n+ 1.$

Our next aim is to establish  a {\em 2nd Main Theorem} for $\exh_A f$  and
its $k$-jet lifts 
 with truncated counting functions
at level one.
\end{abstract}



\section{Introduction}\label{intr}
Let $A$ be a semi-abelian variety of dimension $n\,(>0)$
 with  Lie algebra $\Lie(A)$ and an exponential map
$\exp_A: \Lie(A) \to A$.
Let $f: \C \to \Lie(A)$ be an entire curve
 (i.e., a holomorphic map from $\C$ into $\Lie(A)$),
  and set
\begin{equation}\label{exhat}
\exp_A f:=\exp_A \circ f: \C \to A, \quad
\exh_A f:= ( \exp_A f, f): \C \to A \times \Lie(A).
\end{equation}
We denote by $\overline{\exh_A f (\C)}^{\zar}$  the Zariski closure
  of $\exh_A f (\C)$ in $A \times \Lie(A)$. Put
\[
\td_\C\, \exh_Af:= \dim_\C \overline{\exh_A f (\C)}^{\zar}.
\]
We say that $f$ is  ``{\em $A$-(resp.\ non)degenerate}''  if
\begin{llst}\label{gcond}\rm\begin{itemize}\item[]\hskip-27pt
there is a (resp.\ no) connected  algebraic proper subgroup
 $G \subsetneq A$ such that $G +\exp_Af(0)\supset \exp_A f(\C)$.
\end{itemize}
\end{llst}

\noindent
{\bf(1)} 
It is our first aim to  give a Nevanlinna theoretic proof to
the analytic version of the Ax-Schanuel Theorem for semi-abelian varieties.
\begin{thm}\label{axsch}
If an entire curve $f: z \in \C \to \Lie(A)$ with a semi-abelian variety
 $A$ is $A$-non\-degenerate, 
then 
$\td_\C \, \exh_A f \geq n+ 1.$
\end{thm}

The present study is motivated by S. Lang \cite{la66},
 and J. Ax \cite{ax71}, \cite{ax72}, 
 where Ax proved the above statement for {\em formally analytic}
maps $f: \C \to \Lie(A)$, i.e., $f$ represented by a
vector of {\em formal power series} ({\bf Ax-Schanuel Theorem});
he dealt with the case of formal power series of several variables.
The proof of  Ax \cite{ax71}, \cite{ax72}
 is   based on Kolchin's theory of differential algebra. 
Note that for $f: \C \to \C^n$ and 
$\exp_{(\C^*)^n}: \C^n \ni(w_j) \mapsto (e^{w_j}) \in (\C^*)^n$, $f$ is
 $(\C^*)^n$-nondegenerate
if and only if $f_j(z)-f_j(0) ~(1 \leq j \leq n)$ are $\Q$-linearly
independent.

Because of the nature of the statement, it is natural and
 interesting to
look for an analytic proof in the analytic setting.
The present proof relies on the Log Bloch-Ochiai Theorem \cite{n77}, \cite{n81},
 and Nevanlinna theory of entire curves into semi-abelian varieties
(cf.\ Noguchi-Winkelmann \cite{nw14} in general).

The above Ax-Schanuel Theorem was obtained as an analogue of
the following attractive conjecture:

\begin{llst}\label{schc}\rm {\bf Schanuel's Conjecture}
 (cf.\ Lang \cite{la66} p.\,30). 
Let $\alpha=(\alpha_j) \in \C^n$ be a vector such that
$\alpha_j ~(1\leq j \leq n)$ are $\Q$-linearly independent,
and set $\exh\, \alpha=( (e^{\alpha_j}), \alpha) \in (\C^*)^n \times \C^n$.
Then, 
$ \td_\Q \exh\, \alpha \geq n.$
\end{llst}

This conjecture is known only for $n=1$ (Gel'fond-Schneider
(\cite{la66}, \cite{wal00}); Hilbert's 7th Problem), and
even in $n=2$ it implies the algebraic independence of
$e$ and $\pi$;  cf.,  e.g., Waldschmidt \cite{wal00} \S1,
 Jones-Wilkie \cite{jw15} for more.

\smallskip
\noindent
{\bf(2)}  It is our second aim to establish a 2nd Main Theorem for
 $\exh_Af: \C \to  A \times \Lie(A)$ and its jet lifts.
 Let $J_k(\exh_Af): \C \to J_k(A \times \Lie(A))$ be the $k$-jet lift
of $\exh_Af$ and let $X_k(\exh_Af)$ be the Zariski closure of the image
$J_k(\exh_Af)(\C)$ in $J_k(A \times \Lie(A))$. Because of the flat
 structure of $J_k(A)\iso A \times J_{k,A}$ we may write
\begin{align}
X_k(\exh_Af) \subset A \times \Jk \subset J_k(A \times \Lie(A)),
 \quad \Jk:=\Lie(A)\times J_{k,A}
\iso \C^n \times \C^{nk},
\end{align}
where $J_{k,A}$ is the jet part of $J_k(A)$ (see \S\ref{jb} for details).
Let $\Jkb$ be a projective compactification of $\Jk$; e.g.,
$\Jkb=\pnc \times \P^{nk}(\C)$.

The main result is as follows (cf.\ \S\S\ref{ordf}, \ref{nevt} for notation):
\begin{thm}[2nd Main Theorem]\label{smt}
Let $A$ be a semi-abelian variety.
Let $\exh_A f: \C \to  A \times \Lie(A)$ be the entire curve associated
with an $A$-nondegenerate entire curve $f: \C \to \Lie(A)$.
Then we have:
\begin{enumerate}
\item
Let $Z$ be a reduced algebraic  subset  of $X_k(\exh_A f) ~(\subset A \times \Jk)$
$(k \geqq 0)$.
Then there exists an equivariant projective smooth compactification $\bar A$
of $A$ with the  closure $\bar X_k(\exh_A f)$ (resp.\ $\bar Z$)
 of $X_k(\exh_A f)$ (resp.\ $Z$) in $\bar A \times \Jkb$ such that
\begin{equation}
\label{smtinq}
T_{J_k(\exh_Af)}(r, \omega_{\bar Z}) \leq N_1 (r, J_k(\exh_A f)^* Z)+
 S_{\ep, \exp_Af}(r).
\end{equation}
\item
Moreover, if $\codim_{X_k(\exh_A f)}Z \geqq 2$, then
\begin{equation}\label{cd2}
T_{J_k(\exh_A f)}(r, \omega_{\bar Z}) = S_{\ep, \exp_Af}(r).
\end{equation}
\item\label{k0}
{\rm ($k=0$)} In particular, if $D$ is a reduced effective algebraic divisor on
 $A \times \Lie(A)$ such that $D \not\supset X_0(\exh_A f)$,
then there is an equivariant projective smooth compactification $\bar A$ of
     $A$ such that 
\begin{equation}
T_{\exh_Af}(r, L(\bar D))\leq N_1(r, (\exh_Af)^*D)+
 S_{\ep, \exh_Af}(r),
\end{equation}
where  $\bar D \subset \bar A \times \overline{\Lie(A)}$
with  $\overline{\Lie(A)}:= \Job$.
\end{enumerate}
\end{thm}

As an application we obtain the following (cf.\ 
 \ref{lc} (ii), 
 \ref{bgl} \eqref{cn}
in \S\ref{nevt} below, and
Corvaja-Noguchi \cite{cn} for entire curves into semi-abelian varieties):
\begin{thm}\label{intd}
Let $\exh_Af: \C \to A$ be as in Theorem \ref{smt} and 
let $D$ be an $A$-big divisor on $X_0(\exh_Af)$.
 Then there is an irreducible component $E$ of $D$
such that $\exh_Af(\C) \cap E$ is Zariski dense in $E$;
in particular, the cardinality $|\exh_Af(\C) \cap D|=\infty$.
\end{thm}

Here, $D$ being $A$-big  means roughly that $D$ contains a big divisor
in $A$-factor (see Definition \ref{Abig}).

The present paper is  organized so that  in \S\ref{ordf}
 we prepare the notation from Nevanlinna theory and prove
a  key estimate lemma (see Lemma \ref{keyest}).
We give a proof of  Theorem \ref{axsch} in \S\ref{prf} and those of
Theorems \ref{smt} and \ref{intd} in \S\ref{nevt}.
In \S\ref{exm} we discuss some examples.
In \S\ref{rmk} we remark the generalizations from $\C$ to affine
algebraic curves and other domains ($1$-dimensional).

\begin{rmk}\rm
We would like to refer to some other results  related to the present subject
and remarks.
\begin{enumerate}
\item
J. Tsimerman \cite{ts15} gave another proof to the above
Ax-Schanuel Theorem \cite{ax71} for $A=(\C^*)^n$ by means of
the  ``o-minimal structure'' theory, and similarly Y. Peterzil and
S. Starchenko \cite{ps18} for semi-abelian varieties.
\item
In Noguchi \cite{n18}, there is a direct application of
the Log Bloch-Ochiai Theorem  to the proof
of Raynaud's Theorem (Manin-Mumford Conjecture) through the aid of the
 ``o-minimal structure'' theory.
\item
In the  analytic theory we can think of the ``values'';
 it might be an advantage of the analytic theory, compared with the
     formally analytic theory.
\end{enumerate}
\end{rmk}

{\em Acknowledgment.}  The author is grateful to Professors P. Corvaja
and U. Zannier for interesting discussions on topics related to
Diophantine approximation.

\section{Order functions and a key lemma}\label{ordf}
\subsection{Order functions}\label{ord}
In general for order functions, cf.\ Hayman \cite{h64} Chap's. 1, 2,
Noguchi-Ochiai \cite{no} Chap.\,V, 
and Noguchi-Winkelmann \cite{nw14} Chap.\,2.

Let $X$ be a compact complex space with reduced structure sheaf $\O_X$,
and let $\sI \subset \O_X$ be a coherent ideal sheaf.
Let $f: \C \to X$ be an entire curve. We are going to define an
order function $T_f(r, \omega_\sI)$ of $f$ with respect to $\sI$
(see \,\cite{nw14} \S2.4 for details:
Note that the projective algebraic condition on $X$ is unnecessary because
of $1$-dimensionality of the domain $\C$).

The pull-back $f^*\sI$ is an effective divisor on $\C$, unless
$f(\C) \subset \supp \sI ~(:=\{x \in X: \sI_x \not=\O(X)_x\})$, which we assume.
Denoting by $\ord_z f^*\sI$ the order of $f^*\sI$ at $z \in \C$,
we define the counting functions of $f^*\sI$ truncated at level
 $l \in \N\cup\{\infty\}$  by
\[
 n_l(t, f^*\sI)=\sum_{|z|<t}\min\{ \ord_z f^*\sI, l\}, \quad
 N_l(t, f^*\sI)=\int_1^r \frac{ n_l(t, f^*\sI)}t dt, \quad r>1.
\]
Let $\phi_\sI(x) ~(x \in X)$ be the proximity potential of $\sI$
(see \cite{nw14} \S2.4).
Then we have the pull-back $f^*\phi_\sI(z)$ such that it is of $C^\infty$-class
outside the singular set $f^{-1}(\supp \sI)$ and at a singular point
$a \in f^{-1}(\supp \sI)$ it is written locally as
\[
 f^*\phi_\sI(z)=\lambda_a \log|z-a| + \hbox{$C^\infty$-term},
\]
where $\lambda_a \in \N$ (positive integers).
 We set the proximity function of $f$ for $\sI$ by
\[
 m_f(r, \sI)=\frac1{2\pi}\int_{|z|=r} -f^*\phi_{\sI}(re^{i\theta})
 d\theta.
\]
With $\omega_\sI:=\frac1{\pi i}\del\delbar \phi_{\sI}$ 
 we define the order function of $f$ with respect to $\omega_\sI$ by
\[
T_f(r, \omega_{\sI})=\int_1^r \frac{dt}{t}\int_{|z|<t}
f^*\omega_{\sI}
\quad (r>1).
\]
We have a so-called First Main Theorem for $f$ and $\sI$ (see
\cite{nw14} Theorem 2.4.9):
\begin{equation}\label{fmt}
T_f(r, \omega_\sI)=N_\infty(r,f^*\sI)+m_f(r, \sI)-m_f(1, \sI).
\end{equation}

In the case where $\sI=\sI\ang{Y}$ is the ideal sheaf of an analytic
subset $Y$ possibly non-reduced of $X$, we write
\[
f^*\sI\ang{Y}=f^*Y, \quad \omega_{Y}=\omega_{\sI\ang{Y}}, \quad
 T_f(r, \omega_{Y})=T_f(r, \omega_{\sI\ang{Y}}).
\]

If $X$ is smooth and $Y$ is a divisor $D$ on $X$. Then we have the line
bundle $L(D) \to X$ associated with $D$, and the first Chern class
$c_1(D)$. Then, $\omega_D \in c_1(D)$, and it is common to write
\[
 T_f(r, \omega_D)=T_f(r, L(D))=T_f(r, c_1(D))
\]
for the order functions.

Assume that $X$ is projective algebraic. Let $D$ and $D'$ be  big divisors on
$X$ such that $f(\C) \not\subset \supp D \cup \supp D'$.
Then we have
\[
 T_{f}(r, L(D))=O\left(T_{f}(r, L(D'))\right),
\quad T_{f}(r, L(D'))=O\left(T_{f}(r, L(D))\right).
\]
Therefore in an estimate such as $O\left(T_{f}(r, L(D))\right)$ the
choice of of $D$ or $L(D)$ does not matter; in such a case we simply write
$T_{f}(r)$ for $T_{f}(r, L(D))$ with respect to
some ample or big line bundle $L(D)$ over $X$ as far as
$f(\C) \not\subset \supp D$;
 however, once it is chosen, it is fixed.

Note that the followings are equivalent for $f: \C \to X$:
\begin{llst}\rm\label{rat}\begin{enumerate}\item
 $f$ is rational (not transcendental);
\item
$ \lowlim_{r \to \infty} \frac{T_f(r)}{\log r} < \infty$.
\item
$ \upplim_{r \to \infty} \frac{T_f(r)}{\log r} < \infty$; i.e.,
$T_f(r)=O(\log r)$.
\end{enumerate}
\end{llst}
The order $\rho_f $ of $f: \C \to X$ is defined by
\[
0\leq  \rho_f= \upplim_{r \to \infty} \frac{\log T_f(r)}{\log r} \leq \infty.
\]
If $\rho_f < \infty$, $f$ is said to be of finite order and
\[
 T_f(r)=O(r^{\rho_f + \ep}), \qquad \fa \ep>0.
\]

\subsection{Key lemma}\label{keyl}
 Let $A$ be a semi-abelian variety.
Here we fix an isomorphism $\Lie(A)\iso \C^n$ with coordinates
$(z_1, \ldots, z_n)$.
Let $f=(f_1, \ldots, f_n) : z \in \C \to f(z) \in \C^n (\iso \Lie(A))$
 be an entire curve. In the present paper we assume that $f$ is
 non-constant. 
We denote by $T(r, f_j)$ Nevanlinna's order function 
 (see \cite{h64}, \cite{nw14} Chap.\,1), and  set
\[
 T_f(r)= \max_{1 \leq j \leq n} T(r, f_j).
\]

Let $\bar A$ be an equivariant smooth projective  compactification of $A$
with boundary divisor $\del A=\bar A \setminus A$
of simple normal crossing type (cf.\ \cite{nw14} Chap.\,5),
and let $L \to {\bar A}$ be an ample line bundle.
We set the order function of $\exp_A f: \C \to A$ by
\[
 T_{\exp_Af}(r)=  T_{\exp_Af}(r, L).
\]
For $\exh_A f=(\exp_Af, f): \C \to A \times \C^n$ we set
\[
 T_{\exh_A f}(r):= T_{\exp_Af}(r)+ T_f(r).
\]
By \cite{nw14} Theorem 6.1.9 and \ref{rat} above we see
\begin{prop} Let the notation be as above.
The followings are equivalent:
\begin{enumerate}
\item
$f: \C \to \C^n$ is rational (i.e., $f_j$ are polynomials).
\item
$\exp_A f$ is of  finite order.
\item
$\exh_A f$ is of finite order.
\end{enumerate}
\end{prop}

As usual in Nevanlinna theory we use the symbol
 ``$S_{\exp_Af}(r)$'' for a small term such that
\begin{equation}\label{Sr}
S_{\exp_Af}(r) = O(\log^+ T_{\exp_Af}(r)) +O(\log r)+O(1) ~~ ||,
\end{equation}
where ``$||$'' stands for the validity of the estimate except for $r$
in  exceptional intervals with finite total length, and 
for $\exp_A f$ of finite order there are no such exceptional intervals.
We use the notation through the paper.

The following is the  key lemma for the estimates in the
arguments henceforth.
\begin{lem}\label{keyest}
With the notation  above, we have:
\begin{enumerate}
\item
$ T_f(r) = S_{ \exp_Af}(r)$.
\item
$T_{\exh_Af}(r) = T_{\exp_Af}(r)+S_{\exp_A f}(r) $.
\end{enumerate}
\end{lem}
\begin{pf} It suffices to prove (i).
When the order of $\exp_Af$ is finite, then $f$ is rational, and so
$T_f(r)=O(\log r)$ without exceptional intervals.

In general, we take a representation of the semi-abelian variety $A$
\[
 0 \to (\C^*)^p \to A \to A_0 \to 0,
\]
where $A_0$ is an abelian variety. Then $A$ has a  structure of
locally flat $(\C^*)^p$-principal bundle with transition transformation
by $(S^1)^p:=\{ \zeta \in \C: |\zeta|=1\}^p$-multiplication
(cf.\ \cite{nw14} \S6.1).
After a change of indices of the coordinates $(z_j)$ and
 a linear transform of $(z_j)$ we have the following
expression of the order function $T_{\exp_Af}(r)$ (cf.\ ibid.):
\begin{align}\label{repf}
f(z) &=(f_1(z), \ldots, f_p(z), f_{p+1}, \ldots, f_n(z)), \\
\notag
T_1(r) &:= \sum_{j=1}^p T\left(r, e^{f_j}\right), \\
\notag
T_2(r) &:= \frac1{4\pi}\int_{|z|=r} \sum_{j+1}^n |f_j(z)|^2 \, d\theta
- \frac1{4\pi}\int_{|z|=1} \sum_{j+1}^n |f_j(z)|^2 \, d\theta,\\
\notag
T_{\exp_Af}(r) &=T_1(r) + T_2(r).
\end{align}

For $f_j ~(p+1 \leq j \leq n)$ we have
\begin{align*}
T(r, f_j) &= \frac1{2\pi} \int_{|z|=r} \log^+|f_j| \,d\theta 
\leq  \frac1{4\pi} \int_{|z|=r} \log( 1+ |f_j|^2) \,d\theta \\
&= \frac1{2} \log \left(1+ \frac1{2\pi} \int_{|z|=r} |f_j|^2 \,d\theta
 \right)
\leq \frac12 \log^+T_2(r)+O(1)\\
&= O(\log^+T_{\exp_Af}(r)) +O(1).
\end{align*}

For $f_j  ~(1 \leq j \leq p)$ we have
\[
 T\left(r, e^{f_j}\right) =\frac1{2\pi} \int_{|z|=r}\log^+
\left|e^{f_j(z)}\right|\,d\theta=
\frac1{2\pi} \int_{|z|=r}  \re^+ f_j(z) \,d\theta,
\]
where  $\re^+f_j:=\max\{\re f_j, 0\}$ with the real part $\re f_j$.
With the imaginary part $\im f_j(0)$ and the complex Poisson
 kernel we write
\[
 f_j(z)= \frac{1}{2\pi} \int_{|\zeta|=R} \frac{\zeta+z}{\zeta-z}
\re f_j(\zeta)\, d\theta + \im f_j(0).
\]
For $|z|=r <R$ we get
\[
 |f_j(z)| \leq \frac{R+r}{R-r}\cdot \frac{1}{2\pi}
 \int_{|\zeta|=R} \re^+ f_j(\zeta)\, d\theta + |\im f_j(0)|.
\]
Then we have
\begin{align*}
T(r, f_j) &\leq \frac{1}{2\pi} \int_{|\zeta|=r} \log
(1+|f_j(z)|)\, d\theta \leq
\log\left(1+ \frac1{2\pi} \int_{|z|=r} |f_j(z)|\, d\theta  \right)\\
&\leq \log \left( 1+ \frac{R+r}{R-r}\cdot \frac{1}{2\pi}
 \int_{|\zeta|=R} \re^+ f_j(\zeta)\, d\theta +O(1)\right)\\
&\leq \log \left( 1+ \frac{R+r}{R-r}\cdot T_1(R) +O(1)\right)\\
&\leq \log^+ \left(\frac{R+r}{R-r}\cdot T_1(R)\right)+O(1)
\end{align*}
Now we take $R=r+1/T_1(r)$, so that
\begin{align*}
T(r, f_j) \leq \log^+ \left((2r+1)T_1(r) \cdot
T_1\left(r+\frac1{T_1(r)}\right)\right)+O(1)
\end{align*}
Borel's Lemma (cf.\  Hayman \cite{h64} Lemma 2.4) implies
\[
 T_1\left(r+\frac1{T_1(r)}\right) \leq 2 T_1(r) \quad ||.
\]
Therefore it follows that
\begin{align*}
T(r, f_j) &=O(\log^+T_1(r))+O(\log r)+O(1) \\
&=O(\log^+T_{\exp_Af}(r))+O(\log r)+O(1) \quad ||.
\end{align*}

The proof is completed.
\end{pf}

\section{Proof of  Theorem \ref{axsch}}\label{prf}
By \cite{n81} we see that the Zariski closure
$\overline{\exp_A f(\C)}^\zar$ in $A$ is a translate of a connected algebraic subgroup
(a semi-abelian subvariety)  of $A$.
It follows from $A$-nondegeneracy \ref{gcond} that
$\overline{\exp_A f (\C)}^\zar=A$, so that
\begin{equation}\label{trd2}
\td_\C \, \exp_A f =n.
\end{equation}
Let $\C(A)$ (resp.\ $\C(f)$) be the rational function field of $A$
(resp.\ the field generated by $f_j ~(1 \leq j \leq n)$ over $\C$).
We denote by $\td_{\C(f)} (\exp_A f)^* \C(A)$ the transcendence degree
of the pull-backed field $(\exp_A f)^* \C(A)$ over $\C(f)$.
We prove:
\begin{lem}\label{tdlem}
With the notation  above we have
\begin{equation}\label{1trans}
\td_{\C(f)} \, (\exp_A f)^* \C(A) \geq 1.
\end{equation}
\end{lem}
\begin{pf}
We take a transcendence basis $\{\phi_j\}_{j=1}^n$ of
$\C(A)$ over $\C$ such that $\hat\phi_j:=\phi_j \circ \exp_A f$
are defined as non-constant meromorphic functions,  and
\[
 \hat\phi:=\left(\hat\phi_1, \ldots, \hat\phi_n\right).
\]			 
Assume contrarily that \eqref{1trans} is false; i.e.,
\begin{equation}\label{assp}
\td_{\C(f)} \hat\phi=0 \quad (\hbox{regarded as } \hat\phi=\{\hat\phi_j\}_{j=1}^n).
\end{equation}
Then all $\hat\phi_j$ are algebraic over $\C(f)$.
There are non-zero polynomials $P_j(t)$ in
one variable with coefficients in $\C(f)$
such that
\begin{equation}\label{algr}
 P_j(\hat\phi_j)=0, \quad 1 \leq j \leq n.
\end{equation}
By \cite{nw14} Lemma 2.5.15 we have
\[
 T(r, \hat\phi_j)=O(T_f(r)) +O(1).
\]
With setting $\what T(r):=\max_{1 \leq j \leq n}  T(r, \hat\phi_j)$ we
 thus obtain
\[
 \what T(r) =O(T_f(r))+O(1).
\]
On the other hand, it follows  from \cite{nw14} Theorem 2.5.18 that
\[
 T_{\exp_Af}(r)=O(\what T(r) )+O(1).
\]
Therefore we see that
\[
 T_{\exp_Af}(r)=O( T_f(r) )+O(1).
\]
But this contradicts Lemma \ref{keyest}.
\end{pf}

{\it Continuation of the proof of Theorem \ref{axsch}}:
By \eqref{trd2}, $\td_{\C} \{f, \hat\phi \} \geq n$.
For proof by contradiction we assume that
\[
\td_\C \, \{ f, \hat\phi \} = n.
\]
Then all $f_j$ are algebraic over $\C(\hat\phi)$,
so that there are non-trivial algebraic relations,
\begin{equation}\label{algr1}
P_j(f_j, \hat\phi)=P_j(f_j, \hat\phi_1, \ldots, \hat\phi_n)=0,
 \quad 1 \leq j \leq n.
\end{equation}

If $\td_\C\{f_j\}_{j=1}^n=n$,
the assumption implies  $\td_{\C(f)} \{f, \hat\phi\}=0$;
this does not take place by Lemma \ref{tdlem}.
Therefore $\td_\C\{f_j\}_{j=1}^n < n$, and hence
 there is a non-trivial algebraic relation over $\C$:
\begin{equation}\label{algr2}
Q(f_1, \ldots, f_n)=0.
\end{equation}
If, to say, $f_1$ is contained in \eqref{algr2}, we take the
resultant of \eqref{algr2} and \eqref{algr1} ($j=1$) with respect to $f_1$,
which yields a non-trivial algebraic relation
\[
 Q_1(f_2, \ldots, f_n, \hat\phi)=0.
\]
After repeating this process at most $n$-times we eliminate
$f_1, \ldots, f_n$ in \eqref{algr2} to obtain
a non-trivial algebraic relation
\[
 \hat Q \left( \hat\phi_1, \ldots, \hat\phi_n \right)=0\,;
\]
this again contradicts the algebraic independence of
$\hat\phi_1, \ldots,\hat\phi_n$.

The proof of Theorem \ref{axsch} is completed. \qed

\begin{rmk}\rm\begin{enumerate}
\item
 It is noticed that the logarithmic function $\log (1+t) \in \C[[t]]$
can be dealt with in the Ax-Schanuel Theorem as a formal power series, but cannot
 in our  Theorem \ref{axsch}. To deal with the case of a finite unit disk as
 a domain instead of $\C$ we need some growth condition for
$T_{\exp_A f}(r)$ (cf.\ \S\ref{fdisk}). 
\item
Let $\exp_A: \C^n =\Lie(A) \to A$ be as above.  We have a
     semi-lattice $\Lambda=\Ker \exp_A \subset \C^n$ (the periods of $A$).
Then an entire curve $f: \C \to \C^n$ is $A$-(resp.\ non)degenerate
 if and only if there is a (resp.\ no) complex vector
subspace $E \subsetneq \C^n$ such that $E \supset (f(\C)-f(0))$ and
$E/(E\cap\Lambda)$ is a semi-abelian
variety. Therefore if $\Lambda$ is concerned, it would be better to
say  $f$ being {\em $\Lambda$-(non)degenerate}.
\end{enumerate}
\end{rmk}

\section{Nevanlinna theory of entire curves $\exh_A f$}\label{nevt}
\subsection{Back ground}\label{bg}
  In the same monograph \cite{la66} (1966) as 
 Schanuel's Conjecture \ref{schc} was mentioned, S. Lang raised 
an interesting question (p.\,32):
\begin{llst}\label{lc}\rm\begin{enumerate}
\item\label{lc1}
Let $\vp: \C \to A$ be a $1$-parameter subgroup of an abelian variety $A$
(say Zariski dense), and let $D$ be a hyperplane section of $A$.
Then, is $\vp(\C) \cap D \not= \emptyset$?
\item\label{lc2}
And unless $\psi$ is algebraic, is the cardinality
$|\vp(\C) \cap D| = \infty$?
\end{enumerate}
\end{llst}

It has developed roughly as follows (a non-complete list):
\begin{llst}\label{bgl}\rm
\begin{enumerate}
\def\labelenumi{\rm(\arabic{enumi})}
\item
J. Ax \cite{ax72} (1972) gave an affirmative answer to \ref{lc}\,(i) above.
\item
P.A. Griffiths \cite{gr72} (Problem F, 1972)
 generalized problem \ref{lc}\,(i) for entire curves
into $A$ (so-called Lang's Conjecture).
\item
At the Taniguchi Symposium ``Geometric Function Theory in Katata 1978
organized by S. Murakami (chair), 
 the author formulated a 2nd main theorem for entire curves
 $f: \C \to A$ and a divisor $D$ on $A$ as a conjecture,
 which implies  (2) above and \ref{lc}\,(i) as well; see
 Noguchi--Ochiai \cite{no} (p.\,248, 1984/'90).
\item
Siu and Yeung \cite{sy} (1996) solved 
Lang's Conjecture (2) above 
 for entire curves into abelian  varieties,
and Noguchi \cite{n98} (1998) generalized it for entire curves into semi-abelian
     varieties with another proof, which unifies the result for abelian varieties
 and the classical E. Borel's results for $(\C^*)^n$.
\item\label{smt0}
Noguchi, Winkelmann and Yamanoi \cite{nwy00} (2000), \cite{nwy02} (2002)
proved (3) above, the 2nd Main Theorem for entire curves into semi-abelian varieties,
and finally in \cite{nwy08} proved it with counting functions truncated
     at level one.
\item\label{cn}
P. Corvaja and J. Noguchi \cite{cn} (2012) solved affirmatively \ref{lc}\,(ii) for
entire curves into semi-abelian varieties, $f: \C \to A$ by making
use of the 2nd Main Theorem of \eqref{smt0} above.
It is noticed that \ref{lc}\,(ii) had been open even for $1$-parameter subgroups of
abelian varieties.
\end{enumerate}
\end{llst}

It is natural and interesting to ask  questions similar to the above
 for
\[
 \exh_Af: z \in \C \to (\exp_Af (z), f(z)) \in   A \times \Lie(A)
\]
in view  of
the analytic Ax-Schanuel Theorem \ref{axsch} and
 Schanuel's Conjecture \ref{schc}.

\subsection{2nd Main Theorem}
We denote by $ S_{\ep, \exp_Af}(r) ~(\geq 0)$
 a small term such that
for every $\ep>0$
\[
  S_{\ep, \exp_Af}(r) \leq \ep T_{\exp_Af}(r)+O(\log r) ~~ ||_\ep
\]
(cf.\ \eqref{Sr}).

The main aim is to prove Theorem \ref{smt} for
 $\exh_Af: \C \to A \times \Lie(A)$, which is
 very analogous to \cite{nw14} Theorem 6.5.1
 (cf.\ \cite{nwy08}) for $\exp_Af: \C \to A$.
The proof of Theorem \ref{smt} is in fact an adaptation of the arguments
in \cite{nw14} Chap.\,6 (cf.\ \cite{nwy02}, \cite{nwy08}) by making use
of Lemmata \ref{keyest}, \ref{estIk} for 
the projection $\what I_k$ of \eqref{Ik} below.
Henceforward we will sketch the key points.

The way to obtain the compactifications of $\bar A$ and $X_k(\exh_Af)$
in Theorem \ref{smt} is not written precisely in
 \cite{nw14} Theorem 6.5.1, but it follows from the arguments of the
proof there.

\subsection{Reduction}\label{red}
 Let $f: \C \to \Lie(A)$ be an entire curve.
By the Log Bloch--Ochiai Theorem (\cite{nw14} Theorem 6.2.1,
\cite{n77}, \cite{n81}) $\overline{\exp_A f(\C)}^\zar$ is a translate
of a subgroup $B$ of $A$. By a translate we may assume that
$\overline{\exp_A f(\C)}^\zar=B$. 
Then $f(\C) \subset \Lie(B)~(\subset \Lie(A))$,
and
\[
 \exh_A f: \C \lto B \times \Lie(B) \subset A \times \Lie(A).
\]
Now, $f: \C \to \Lie(B)$ is $B$-nondegenerate.
Therefore without loss of generality
 we may assume  that $B=A$, i.e, $f$ is $A$-nondegenerate.

\subsection{Jet bundles}\label{jb}
We keep the same notation as in the previous subsections.
Let $f: \C \to \C^n$ be an $A$-nondegenerate entire curve.
We would like to study the value distribution of
 $\exh_A f: \C \to A \times  \C^n$.

Let $J_k(A) \to A$ (resp.\ $J_k(\Lie(A)) \to \Lie(A)$) be
the $k$-th jet bundle over $A$ (resp.\ $\Lie(A)$).
Because of the flat structure of the logarithmic tangent
(and cotangent, as well) bundle over $A$ (cf.\ \cite{nw14} \S4.6.3),
 we have the trivializations:
\begin{align}
\begin{array}{cc}
 J_k(A) \iso A \times J_{k,A}, & J_{k, A} \iso \C^{nk},\\
 J_k(\Lie(A)) \iso \Lie(A) \times J_{k,\Lie(A)}, &
 J_{k,\Lie(A)} \iso \C^{nk},
\end{array}
\end{align}
where $J_{k, A}$ and $J_{k,\Lie(A)}$ are the so-called
jet-parts of $A$ and $\Lie(A)$, respectively.
Through the exponential map $\exp_A: \Lie(A) \to A$ we have the
natural isomorphism $J_{k, A} \iso J_{k,\Lie(A)}$, which
are identified. Therefore we have
\begin{equation}
J_k(A \times \Lie(A)) \iso A \times \Lie(A) \times J_{k,A}\times
J_{k,A}.
\end{equation}
Let $\Delta_k \subset  J_{k,A}\times J_{k,A}$ be the diagonal,
and let $J_k(\exh_Af): \C \to J_k(A \times \Lie(A))$ be the $k$-jet
lift of $\exh_Af$. Then we see that
\begin{equation}
J_k(\exh_A  f): z \in \C \to (\exp_A f(z), f(z), J_{k, f}(z), J_{k,f}(z))
\in A \times \Lie(A) \times \Delta_k,
\end{equation}
where $J_{k, f}$ is the jet part of $J_k(f)=(f, J_{k, f}) \in \Lie(A)
\times J_{k, A}$. For the sake of simplicity we identify
$\Delta_k=J_{k,A}$ and write
\begin{equation}
J_k(\exh_A  f): z \in \C \to (\exp_A f(z), f(z), J_{k, f}(z))
\in A \times \Lie(A) \times J_{k, A} (\subset J_k(A \times \Lie(A))).
\end{equation}
We put
\begin{align}
\Jk &= \Lie(A) \times J_{k,A} ~ \iso \C^n \times \C^{nk} \quad
(\hbox{\em extended jet part}), \\
 X_k(\exh_Af) &=\overline{J_k(\exh_Af)(\C)}^{\mathrm{Zar}}
\subset A \times \Jk .
\notag
\end{align}
We define the {\em extended jet projection} by
\begin{equation}\label{Ik}
 \what I_k : X_k(\exh_Af) ~(\subset A \times \Jk) \lto \Jk,
\end{equation}
which will play the role of the jet projection
$I_k$ (cf.\ \cite{nw14} p.\,151) for  entire curves into semi-abelian varieties
in \cite{nw14} Chap's.\,4--6.

Let $L_{\bar A} \to \bar A$ be a big line bundle over a
projective compactification $\bar A$ of $A$.
We take a compactification $\Jkb$ of $\Jk$, e.g.,
\[
 \Jkb=\overline{\C^n \times \C^{nk}}=\pnc \times \P^{nk}(\C),
\]
and the ample line bundle $H=O_\pnc (1) \tensor O_{\P^{nk}(\C)}(1) \to \Jkb$,
with which we define
\[
  T_{\what I_k\circ J_k(\exh_Af)}(r) = T_{\what I_k\circ J_k(\exh_Af)}(r, H).
\] 
\begin{lem}\label{estIk}
For $\what I_k$  we have
\begin{equation}\label{estIkeq}
 T_{\what I_k\circ J_k(\exh_Af)}(r) =S_{\exp_Af}(r).
\end{equation}
\end{lem}
\begin{pf}
Since
\[
\what I_k\circ J_k(\exh_Af): z \in \C \to (f(z), J_{k,\exp_Af}(z)) \in
\Lie(A) \times J_{k,A} ,
\]
it follows from  Lemma on logarithmic derivative for $\exp_Af$
(\cite{n77}, \cite{nw14} \S4.7) that
$$ T_{J_{k,\exp_Af}}(r)=S_{\exp_Af}(r).
$$
This combined with Lemma \ref{keyest} implies \eqref{estIkeq}.
\end{pf}

\subsection{$A$-action}
 We consider an $A$-action on
 $A \times \Lie(A)\times J_{k,A} \subset J_k(A \times \Lie(A))$ by
\[
 (a, (x, v, w )) \in A \times \left(A \times \Lie(A)\times J_{k,A}\right)
\to (a+x, v,w) \in A \times \Lie(A)\times J_{k,A}.
\]
We denote the stabilizer subgroup of $X_k(\exh_Af)$ by
\[
\st(X_k(\exh_Af))=\st_A(X_k(\exh_Af))=\{a \in A: a+X_k(\exh_Af)=X_k(\exh_Af)\},
\]
and by $\st(X_k(\exh_Af))^0$ the identity component.

\begin{lem}\label{stX}
With the notation above, $\st(X_k(\exh_Af))^0 \not= \{0\}$.
\end{lem}
\begin{pf}
If otherwise, $\st(X_k(\exh_Af))^0 = \{0\}$.
We consider the $l$-jet space $J_l(X_k(\exh_A f))$ of $X_k(\exh_Af)$
(``jet of jet'') with induced projection
\[
 d^l\what I_k: J_l(X_k(\exh_A f)) ~(\subset J_l(A \times \Lie(A)\times J_{k,A}))
\to J_l(\Lie(A)\times J_{k,A}).
\]
By \cite{nw14} Lemma 6.2.4, there is a large number $l \in \N$ such that
the differential $d(d^l \what I_k)$ is non-degenerate at general points
of $J_l(X_k(\exh_A f))$. Therefore we have
\begin{align*}
T_{\exp_Af}(r) & \leq T_{\exh_Af}(r)=O\left(T_{J_l(J_k(\exh_{A} f))}(r)\right)
 =O\left(T_{d^l\what I_k(J_{k}(\exh_A f))}(r) \right).
\end{align*}
On the other hand,
 $T_{d^l\what I_k(J_{k}(\exh_A f))}(r) =S_{\ep, \exp_Af}(r)$
by Lemma \ref{estIk}; it is a contradiction.
\end{pf}

\begin{prop}[cf.\ \cite{nw14} Theorem 6.2.6]
Let $B=\st(X_k(\exh_Af))^0$ and set the quotient map
\begin{equation}\label{qb}
q_B: X_k(\exh_A f) \to X_k(\exh_Af)/B  
~(\subset (A/B)\times \Lie(A) \times J_{k,A} ~).
\end{equation}
Then $T_{q_B \circ J_k(\exh_A f)}(r)=S_{\exp_Af}(r)$.
\end{prop}
\begin{pf}
The semi-abelian variety $A/B$ acts on 
$(A/B)\times \Lie(A) \times J_{k,A}$ by
 the translations of the first factor and the identity
for the other factors. Then $\st_{A/B}(X_k(\exh_Af)/B)^0=\{0\}$.
As in the proof of Lemma \ref{stX}, with a large $l$ the projection
\[
 \rho_l: J_l(X_k(\exh_A f)/B) \to J_l(\Lie(A) \times J_{k,A})
\]
has a non-degenerate differential $d\rho_l$ at general points
(see \cite{nw14} Lemma 6.2.4).
Therefore we have
\[
 T_{d^lq_B \circ J_k(\exh_Af)}(r)=
O\left(T_{\rho_l \circ d^l q_B \circ J_k(\exh_Af)}(r)\right),
\]
where $d^lq_B: J_l(X_k(\exh_A f)) \to J_l(X_k(\exh_Af)/B) $
is the induced morphism from $q_B$.
It follows from $\rho_l \circ d^lq_B=d^l \what I_k$ and Lemma \ref{estIk} that
\[
 T_{d^lq_B \circ J_k(\exh_Af)}(r)=S_{\exp_A f}(r).
\]
On the other hand we have by Lemma \ref{keyest}
\[
 T_{d^lq_B \circ J_k(\exh_Af)}(r)=T_{q_B \circ J_k(\exh_Af)}(r)
+S_{\exp_Af}(r).
\]
Thus we deduce that $T_{q_B \circ J_k(\exh_Af)}(r)=S_{\exp_Af}(r)$.
\end{pf}

\subsection{Proof of Theorem \ref{smt}}
 Let the notation be as in Theorem \ref{smt}.
Through the arguments of the proof of \cite{nw14} Theorem 6.5.6
with replacing the jet projection $I_k$ there by the extended jet
projection $\what I_k$ (cf.\ \eqref{Ik}) and 
$q_k^B$  by $q_B$ (cf.\ \eqref{qb}), we deduce that
there are a number $l_0 \in \N$ and a compactification $\bar A$ of $A$
 such that
\begin{equation}\label{smtl}
T_{J_k(\exh_Af)}(r, \omega_{\bar{Z}})=N_{l_0}(r, J_k(\exh_Af)^* Z) +
S_{ \exp_Af}(r),
\end{equation}
where $\bar X_k(\exh_Af)$ (resp.\ $\bar Z$)  is the closure of $X_k(\exh_Af)$
(resp.\ $Z$) in $\bar A \times \Jkb$.

Next, we show (ii). It follows from \eqref{smtl} that
\begin{equation*}
T_{J_k(\exh_Af)}(r, \omega_{\bar{Z}}) \leq l_0 N_1(r, J_k(\exh_Af)^* Z) +
S_{ \exp_Af}(r).
\end{equation*}
Now, making use of the assumption $\codim_{X_k(\exh_A f)}Z \geqq 2$
together with  Lemmata \ref{keyest} and \ref{estIk},
we  deduce from the arguments of \cite{nw14}
\S6.5.3, adapted for $\exh_Af: \C \to A \times \Lie(A)$ that
\begin{equation}\label{cd2est}
 N_1(r, J_k(\exh_Af)^* Z) =  S_{\ep, \exp_Af}(r).
\end{equation}
Thus, \eqref{cd2} is deduced, and (ii) is finished.

Now, we go back to the proof of (i).
It follows from  the First Main Theorem \ref{fmt} that
\begin{equation}\label{cd2.1}
 N_1(r, J_k(\exh_Af)^*Z) \leq  N_\infty(r, J_k(\exh_Af)^*Z)\leq
T_{J_k(\exh_Af)}(r, \bar Z). 
\end{equation}
Note that (i) is finished in the case of
 $\codim_{X_k(\exh_Af)} Z \geq 2$.
Thus we consider the case where $Z$ is an effective reduced
 divisor $D$ on $X_k(\exh_Af)$. 
Let $D=\sum D_i$ be the irreducible  decomposition.
We deduce from \eqref{smtl} that
\begin{align*}
 T_{J_k(\exh_Af)}(r, \omega_{\bar D})&\leq N_{l_0}(r, J_k(\exh_Af)^*D)+
S_{\exp_Af}(r)\\
& \leq N_{1}(r, J_k(\exh_Af)^*D)+ 
l_0 \sum_{i<j} N_{1}(r, J_k(\exh_Af)^* (D_i \cap D_j))\\
&\quad + 
l_0 \sum_i N_{1}(r, J_{k+1}(\exh_Af)^* J_1(D_i)) +
S_{ \exp_Af}(r)
\end{align*}
(cf.\ \cite{nw14} (6.5.51)).
Since $\codim_{X_k(\exh_Af)} D_i \cap D_j \geq 2$, \eqref{cd2est} implies
\[
  N_{1}(r, J_k(\exh_Af)^* (D_i \cap D_j))= S_{\ep, \exp_Af}(r).
\]

We have 
$J_{k+1}(\exh_Af): \C \to X_{k+1}(\exh_Af)\subset J_{k+1}(A \times \Lie(A))$
 and $B=\st( X_{k+1}(\exh_Af))^0$.
For each $D_i$ we have two cases: (1) $B \subset \st(D_i)^0$ and
(2) $B \not\subset \st(D_i)^0$.
In the first case (1) we have by using $q_B$ that
\[
N_{1}(r, J_{k+1}(\exh_Af)^* J_1(D))=S_{\exp_Af}(r)
\]
(see \cite{nw14} \S6.5.4 (b)).
In the second case (2), we have by \cite{nw14} Lemma 6.5.50
\[
 \codim_{X_{k+1}(\exh_Af)} (X_{k+1}(\exh_Af) \cap J_1(D_i)) \geq 2.
\]
Then it is deduced from \eqref{cd2est} with $k+1$ that
\[
 N_{1}(r, J_{k+1}(\exh_Af)^* J_1(D))= S_{\ep, \exp_Af}(r).
\]
Thus, \eqref{smtinq} follows.

(iii) The case of  $k=0$ is a special case of (i);
 the proof of Theorem \ref{smt} is completed. \qed

\medskip
We consider the fundamental case where $k=0$
 and $Z$ is a
reduced divisor $D$ on $X_0(f)$.
Let $p_1: X_0(f) \,(\subset A \times \Lie(A)) \to A$ be the projection
to the first factor $A$.
\begin{defn}\label{Abig}\rm
We say that $D$ is {\em $A$-big} if for a big divisor $E$ on $A$
(i.e., the closure $\bar E$ in a compactification $\bar A$ of
$A$ is big)
 the complete linear
system $|mD - p_1^*E |$ with large $m \in \N$
 contains an effective divisor on $X_0(f)$.
\end{defn}

If $D$ is $A$-big, then
\begin{align}\label{a-big}
 T_{\exp_A f}(r) &=O(T_{\exh_Af}(r, \omega_{\bar D})) ,\\
 T_{\exh_A f}(r) &=O(T_{\exh_Af}(r, \omega_{\bar D})) \quad ||.
\notag
\end{align}

\begin{cor}\label{nobig}
Let $f: \C \to \Lie(A)$ be an $A$-nondegenerate entire curve and
let $D$ be a reduced divisor on $X_0(f)$.
If $\ord_z (\exh_Af)^*D \geq 2$ for all $z \in \supp (\exh_Af)^*D$
except for finitely points of $\supp (\exh_Af)^*D$, then $D$ is not $A$-big.
\end{cor}
\begin{pf}
If $D$ is $A$-big,
it follows from the 2nd Main Theorem \ref{smt} and \eqref{a-big} that
\begin{align*}
 T_{\exh_Af}(r, \omega_{\bar D}) &\leq N_1(r, (\exh_Af)^*D)+ S_{\ep, \exh_Af}\\
          &\leq \frac12 N_\infty(r, (\exh_Af)^*D)+ S_{\ep, \exh_Af}\\
          &\leq \frac12 T_{\exh_Af}(r, \omega_{\bar D})+ 
     \ep T_{\exh_Af}(r, \omega_{\bar D})+O(\log r) ~~||_\ep.
\end{align*}
This implies a contradiction `$1 \leq \frac12$'.
\end{pf}

\begin{rmk}\rm
This corollary is motivated through a discussion with Corvaja
and Zannier on their related
or analogous results in rational function fields of
Corvaja and Zannier \cite{cz} (see \cite{n97} too).
\end{rmk}

\subsection{Proof of Theorem \ref{intd}}

With the notation of the theorem we set
\[
 Z =\overline{\exh_Af(\C) \cap D}^\zar \subset  D.
\]
If $\codim_{X_0(\exh_Af)} Z \geq 2$, Theorem \ref{smt} (ii) would imply
\[
 N_1(r, (\exh_Af)^* D)= N_1(r, (\exh_Af)^*Z)= S_{\ep, \exp_Af}(r).
\]
It follows from Theorem \ref{smt} (i) and \eqref{a-big} that for every $\ep>0$
\[
T_{\exh_Af}(r, \omega_{\bar D}) \leq  \ep T_{\exh_Af}(r, \omega_{\bar D})+
O(\log r) ~~ ||_\ep:
\]
This is a contradiction.

Therefore, $Z$ has an irreducible component of codimension one in
$X_0(\exh_Af)$, which is an irreducible component of $D$.  \qed

\section{Examples}\label{exm}
To discuss examples it is convenient to write $\Lie(A) \times A$
for $A \times \Lie(A)$, so that in this section we use the notation
\[
 \exh_A f=(f, \exp_Af): \C \lto \Lie(A) \times A;
\]
there will be no confusion.

(a) The optimality of \eqref{trd2}:
  Let $A=(\C^*)^n$ and let $\alpha_j, 1 \leq j \leq n$,
be complex numbers, linearly independent over $\Q$.
Then the entire curve
$\phi(z)=(\alpha_1 z, \ldots , \alpha_n z) ~(\in \C^n=\Lie(A))$
is $A$-nondegenerate with  the natural exponential map
$\exp :\C^n \ni (z_j) \mapsto (e^{z_j}) \in A$, and so
\[
 \td_\C \{z, e^{\alpha_1z}, \ldots e^{\alpha_n z}\}=n+1.
\]

Let $\pnc \supset \C^n$ and 
 $\bar A:=\ponec^n \supset A$ be the compactifications and let
$ T_{\exp \phi}(r)$ and $ T_{\exh \phi}(r)$ denote
 the order functions with respect to the products
of point bundles. We write $\exh\, \phi=(\phi, \exp \phi)$.
Then,
\[
T_{\exh\, \phi}(r)= T_{\exp \phi}(r) +O(\log r)=
\frac{\sum_{j=1}^n |\alpha_j|}\pi r +O(\log r).
\]
Let $P(z_1, \ldots, z_n, w_1, \ldots , w_n)$ be a polynomial of
  degree $d_j$ in $w_j ~(1 \leq j \leq n)$.
We assume the condition:
\begin{llst}\label{pcond}\rm\begin{enumerate}\item 
 The divisor on $\pnc \times \bar A$ defined by
the zero of $P$ is reduced and equal to the closure 
 $\bar D_P$ of the divisor $D_P$
 defined by $\{P=0\} \cap (\Lie(A) \times A)$.
\item
$D_P$ is $A$-big (see Definition \ref{Abig}).
\end{enumerate}
\end{llst}
\noindent
We need this condition; 
 otherwise, to say, if $P=z_1$, then $D_P=\{z_1=0\}$ and
$D_P$ is not $A$-big; there is only one root $z=0$ in
 $P(\exh\,\phi (z))=0$.
We  have
\[
 T_{\exh\, \phi}(r, L(\bar D_P))= \frac{\sum_{j=1}^n d_j |\alpha_j|}\pi r +O(\log r).
\]
It is yet, in general, hard to find a root
of $P(\exh\,\phi (z))=0$, but by Theorem \ref{smt} (iii)
\begin{align}\label{ordest}
T_{\exh\, \phi}(r, L(\bar D_P)) 
= N_\infty(r, (\exh\,\phi)^* D_P)+  O(\log r)  =N_1(r, (\exh\,\phi)^* D_P)+
  S_{\ep, \exp \phi}(r). 
\end{align}
The above estimate of $N_\infty(r, (\exh\, \phi)^*D_P)$ is classical due to
 Borel-Nevanlinna, but that of $N_1(r, (\exh\, \phi)^*D_P)$ is new;
moreover from
 Theorem \ref{intd} we obtain 
\begin{equation}\label{zdense}
\overline{ \exh \phi(\C) \cap D_P}^\zar=D_P.
\end{equation}

By Corollary \ref{nobig} there is no entire function $g(z)$  such that
\begin{equation}
 P(\exh_Af (z))=(g(z))^m \quad (m \in \N\, \geq 2).
\end{equation}

In view of the transcendence problem
of $\pi$ and $e$ the above example in the case of
$n=2$ and complex vector $\varpi_0=(1, 2\pi i) \in \C^2$
is of a special interest.
We consider the induced $1$-parameter subgroup
\begin{align}\label{spp}
\phi_0(z) &= z \varpi_0 \in \C^2, \qquad z \in \C,\\
 \exh\, \phi_0(z) &=(z, 2\pi i z, e^z, e^{2\pi i z})=
(z_1,z_2,w_1,w_2)\in \C^2 \times (\C^*)^2.
\notag
\end{align}
Let $P(z_1,z_2,w_1,w_2)$ be a polynomial with integral coefficients
satisfying condition \ref{pcond}.
If $P(\exh\, \phi_0(\zeta))=0$, then
$$\zeta, ~ 2\pi i \zeta, ~ e^\zeta,  ~ e^{2\pi i \zeta}$$
 are algebraically dependent, and there are infinitely many such points,
for which \eqref{zdense} and \eqref{ordest} hold.

(c) Let
$\exp_{(j)}(z) ~(j=1,2, \ldots)$ denote the $j$-times iteration
 of the exponential function $e^z$.
We set
\[
f_1(z)=z, ~ f_j(z)=\exp_{(j-1)}(z), ~ 2 \leq j \leq n .
\]
Then we have
\[
 \td_\C \{f_1, \ldots, f_n, e^{f_1}, \ldots, e^{f_n}\} = \td_\C 
\{f_1, \ldots, f_n, f_{n+1}\} = n+ 1.
\]

In this case, $\exp f$ is of infinite order and 
 $T_{\exp f}(r)$ has a  growth such that  $T_{\exp f}(r)\sim \exp_{(n)} (r)$.

(d) (Cf.\ Brownawell-Kubota \cite{bk77})
Let $\exp_A: \Lie(A)\iso \C^n \to A$ be an exponential map of
a semi-abelian variety $A$.
In general if $f_j ~(1 \leq j \leq n)$ are entire functions,
linearly independent over $\C$, then $f=(f_j): \C \to \C^n$
is $A$-nondegenerate. In particular, let $n=l+m$ and
let $\wp_j(w) ~(1 \leq j \leq m)$ be Weierstrass's {\em pe}-functions.
Then
\[
 \td_\C \{f_1, \ldots, f_{l+m}, e^{f_1}, \ldots, e^{f_l},
\wp_1(f_{l+1}), \ldots, \wp_m(f_{l+m})\} \geq l+m+1.
\]

(e) Let $f_1(z)=z, f_2(z)=z$. Then they are not linearly independent
over $\C$. Let $E_j ~(j=1,2)$ be elliptic curves which are not
isogenous to each other. Let $A=E_1 \times E_2$ and
 $\exp_A: \C^2 \to A$ be an exponential map. Then $f=(f_1,f_2): \C \to \C^2$
is $A$-nondegenerate, and so
\[
 \td_\C\{z, \wp_1(z), \wp_2(z)\}=3.
\]
This should be known classically.

Let $\overline{\Lie(E_1)\times\Lie(E_2)}\times A=
\ptwoc \times E_1\times E_2$ be the compactification with  product
line bundle $L$ (resp.\ $L_0$) of the hyperplane bundle and the
 point-bundles over $\ptwoc \times E_1\times E_2$ (resp.\ $E_1 \times E_2$).
 Then, $L$ is ample and
the order functions  satisfy
\[
 T_{\exh f}(r, L)=T_{\exp f}(r, L_0)+ 2\log r +O(1)=\frac{\pi r^2}2
 \left(\frac1{\lambda_1}+\frac1{\lambda_2}+o(1) \right) ,
\]
where $\lambda_j$ is the surface area of the fundamental domain
 of $\wp_j ~(j=1,2)$.
Let $P(z_1, z_2, w_1, w_2)$ be an irreducible
 polynomial, involving $w_1$ and $w_2$, and of
degree $d_1$ (resp.\ $d_2$) with respect to $w_1$ (resp.\ $w_2$).
We consider $w_j=\wp_j ~(j=1,2)$ a rational function of $E_j$, and
denote by $D_P$ the divisor defined by the zeros of $P$ on
 $Lie(A)\times A$. 
Let $\Xi_P$ denote the zero divisor on $\C$ defined by
 $P(z, z, \wp_1(z), \wp_2(z))=0$.
Then we have
\[
N_\infty(r, \Xi_P)=N_1(r, \Xi_P)+  S_{\ep, \exp f}(r)
={\pi r^2}
 \left(\frac{d_1}{\lambda_1}+\frac{d_2}{\lambda_2}+o(1) \right)
+  S_{\exp f}(r). 
\]

It also follows from Theorem \ref{intd} that
$\overline{ \exh_Af (\C) \cap D_P}^\zar$ contains
an irreducible component of $D_P$.
By Corollary \ref{nobig} we see that there is no meromorphic
function $g(z)$ on $\C$ satisfying
\begin{equation}
 P(\exh_A f(z))= g(z)^m \quad (m \in \N,\, \geq 2).
\end{equation}

(f) Set $f_1=z, f_2= z^2, f_3=z$. Then these are not linearly
independent over $\C$. Let $\exp_E: \C \to E$ be an exponential
map of an elliptic curve $E$ with Weierstrass' $\wp(w)$.
Set $A=(\C^*)^2 \times E$ with $\exp_A: \C^3 \to A$.
Then $f=(f_j):\C \to \C^3$  is $A$-nondegenerate, and so
\[
 \td_\C\{z, e^z, e^{z^2}, \wp(z)\}=4.
\]
The order function of $\exh_Af$ has a  growth, 
 $T_{\exh_A f}(r) \sim r^2$.

\section{Remarks to affine algebraic curves and other domains}\label{rmk}

\subsection{Affine algebraic curves}
In this section $A$ denotes a semi-abelian variety.

Let $R$ be a complex affine algebraic curve and let $f: R \to \Lie(A)$ be
a holomorphic curve.
We may consider $\exp_Af: R \to A$ and
 $\exh_Af: R \to A \times \Lie(A)$. 

In general for a holomorphic curve $g: R \to A$ the arguments up to obtaining
an estimate such as \eqref{smtl} work (cf.\ \cite{nw14} Chap.\ 6),
but further to advance to the estimates of the
 2nd Main Theorem \ref{smt} with counting functions truncated at level
 one, we need to lift $g$ to $\tilde g: R \to \Lie(A)$,
which does not exists in general, since $R$ is not simply connected.
But in the present case we begin with a holomorphic curve
 $f:R \to \Lie(A)$, which is a lift of $\exp_Af:R \to A$.
Therefore we can advance the arguments further there.

Let $\bar R$ be the compactification of $R$ by adding a finite number of
points. To study $\exh_A f$  we may localize the problem about 
an infinite point $a \in \bar R\setminus R$; we take a disk neighborhood
$\Delta$ of $a$ in $\Bar R$. Then $\Delta^*=\Delta\setminus\{a\}$
is a punctured disk, and the analysis of transcendental properties of $\exh_Af$
is reduced to that of the restriction $\exh_Af|_{\Delta^*}$
(see the next).

\subsection{Punctured disk}
Let $\Delta^*$ be a punctured disk.
As mentioned above, although $\Delta^*$ is not simply connected,
it does not cause a difficulty here, since we give in first a holomorphic curve
\begin{equation}
f: z \in \Delta^* \lto (f_1(z), f_2(z), \ldots, f_n(z)) \in \C^n
\iso \Lie(A).
\end{equation}
For a notational convenience we put the puncture at infinity and
introduce a coordinate $z$ such that
\[
 \Delta^*=\{|z|>1\}, \qquad \Delta=\Delta^* \cup \{ \infty\}.
\]
Let $F(z)$ be one of $f_j(z)$.
Then $F(z)$ is expanded to a Laurent series
\begin{equation}
F(z) =\sum_{\nu > 0} c_\nu z^\nu + \sum_{\nu \leq 0} c_\nu z^\nu
=F_{\mathrm{m}}(z) + F_{\mathrm{b}}(z).
\end{equation}
Fix $r_0 >1$. Then $F_{\mathrm{b}}(z)$ and their derivatives
are bounded in $\{|z|\geq r_0\}$, and $F_{\mathrm{m}}(z)$ is the
main part of the expansion:
\begin{equation}
\frac{d^k}{dz^k} F(z)=\frac{d^k}{dz^k}F_{\mathrm{m}}(z)+O(1), \qquad
 |z| \geq r_0, ~ k \geq 0.
\end{equation}

Note that $F_{\mathrm{m}}(z)$ is holomorphic in $\C$.
Applying the key Lemma \ref{keyest} for $F_{\mathrm{m}}(z)$, we deduce
the key Lemma \ref{keyest} for $f$ , $\exp_A f$ and $r>r_0$.
We can then deduce the 2nd Main Theorem \ref{smt} for
$\exh_Af: \Delta^* \to A \times \Lie(A)$ and $J_k(\exh_A f)$ as well
for $r >r_0$.

\subsection{Finite disk}\label{fdisk}
We consider the case of a disk of $\C$ with finite radius, to say,
the unit disk $\Delta$.
Let $f: \Delta \to \Lie(A)$ be a holomorphic curve.
In this hyperbolic case, to make the proofs of key Lemma \ref{keyest}
and Lemma on logarithmic derivatives to work at least for a sequence
$r_\nu \nearrow 1 ~(\nu \to \infty)$,
 we need a technical condition
on the growth of the order function $T_{\exp_Af}(r)$ such that
\begin{equation}\label{gwcond}
 \upplim_{r \to 1} \frac{T_{\exp_Af}(r)}{\log \frac1{1-r}}=\infty
\end{equation}
(cf.\ Nevanlinna \cite{ne29} Chap.\ VI, Hayman \cite{h64} \S2.3).
Under this condition the 2nd Main Theorem \ref{smt} for
$\exh_Af: \Delta \to A \times \Lie(A)$ is deduced.

\subsection{Open Riemann surfaces}
Let $R$ be an open Riemann surface.
The generalizations of Nevanlinna theory for meromorphic functions
on $R$,  holomorphic maps from $R$ to another Riemann surface and
holomorphic curves from $R$ into $\pnc$ are
classical (cf., e.g., 
Sario-Noshiro \cite{sn66},
Wu \cite{wu70}).
There one uses  a finite (hyperbolic case) or infinite (parabolic case)
 exhaustion function  $\tau: R \to [0,r_0)$ with $r_0 \leq \infty$ 
such that $\tau$ is harmonic outside a compact subset of $R$.

Similarly it is formally possible to extend the 2nd Main Theorem \ref{smt}
for  holomorphic curves $f:R \to \Lie(A)$ and $\exh_Af: R \to
A\times\Lie(A)$.
There we use the differential $\del \tau$,  holomorphic where
$\tau$ is harmonic. For a logarithmic $1$-form $\omega$ on $A$
we take the ratio $F=(\exp_Af)^*\omega/ \del\tau$,
which may have  poles at zeros of $\del \tau$.
The counting functions of those zeros is the counting function
$E_R(r)$ of the Euler numbers of $\{ \tau < r \}$, which appears
in the estimates \eqref{smtl} and \eqref{smtinq} under a growth assumption such
as  \eqref{gwcond} with $r_0=1$ in hyperbolic case (i.e., $r_0 < \infty$):
\begin{align}
&T_{J_k(\exh_Af)}(r, \omega_{\bar{Z}}) =N_{l_0}(r, J_k(\exh_Af)^* Z) +
C(k,l_0) E_R(r) \\
\notag
&\hskip120pt  + O(\log^+ T_{\exp_Af}(r))+O(\log r) ~~ ||, \\
&(1-\ep)T_{J_k(\exh_Af)}(r, \omega_{\bar Z , J_k(\exh_A f)})
 \leq N_1 (r, J_k(\exh_A f)^* Z) +C'(\ep, k)E_R(r) + O(\log r)~~ ||_\ep.
\end{align}
Here  $C(k,l_0)$ and $C'(\ep, k)$ are
positive constants such that $C(k,l_0), C'(\ep,k) \nearrow \infty$ as
$k, l_0  \nearrow \infty$ and $\ep \searrow 0$; there are no estimates for
 $C(k,l_0), C'(\ep,k)$.
Therefore, in order to obtain a meaningful consequence
 we need  a technical condition 
(besides  \eqref{gwcond} in hyperbolic case) such that
$$\upplim_{r \to r_0} \frac{E_R(r)}{T_{\exp_Af}(r)} =0.$$

\bigskip
\rightline{Junjiro Noguchi}
\rightline{Graduate School of Mathematical Sciences}
\rightline{The University of Tokyo}
\rightline{Komaba, Meguro-ku, Tokyo 153-8914, Japan}
\rightline{e-mail: noguchi@g.ecc.u-tokyo.ac.jp}
\end{document}